\documentclass[a4paper,12pt]{article}

\usepackage{listings}
\usepackage{longtable}

\usepackage{graphicx}

\usepackage[font=scriptsize]{caption}

\usepackage{float}

\usepackage{amsmath,amssymb,mathrsfs,eucal}

\usepackage[T1,T2A]{fontenc}

\font\cyrfont=wncyss10

\def\sza{\hbox{\cyrfont X}} 

\newtheorem{thm}{Theorem}

\newtheorem{conj}[thm]{Conjecture}

\begin{document}

\title{Orders of Tate-Shafarevich groups for the 
cubic twists of $X_0(27)$}

\author{Andrzej D\k{a}browski and Lucjan Szymaszkiewicz} 

\date{}

\maketitle{}

{\it Abstract}. 
This paper continues the authors previous investigations concerning 
orders of Tate-Shafarevich groups in quadratic twists of a given elliptic  
curve, and for the family of the Neumann-Setzer type elliptic curves.  
Here we present the results of our search for the 
(analytic) orders of Tate-Shafarevich groups for the cubic twists of $X_0(27)$. 
Our calculations extend those given by Zagier and Kramarz 
\cite{ZK} and by Watkins \cite{Wat}. Our main observations concern the asymptotic 
formula for the frequency of orders of Tate-Shafarevich groups. 
In the last section we propose a similar asymptotic formula for the class numbers 
of real quadratic fields. 

\bigskip 
Key words: elliptic curves, cubic twists, Tate-Shafarevich group, 
Cohen-Lenstra heuristics, distribution of central $L$-values, class numbers of real 
quadratic fields

\bigskip 
2010 Mathematics Subject Classification: 11G05, 11G40, 11Y50

\section{Introduction}

Let $E$ be an elliptic curve defined over $\Bbb Q$ of conductor $N_E$, 
and let $L(E,s)$ denote its $L$-series. 
Let $\sza(E)$ be the Tate-Shafarevich group of $E$, 
$E(\Bbb Q)$ the group of rational points, and $R(E)$ the regulator,  
with respect to the N\'eron-Tate height pairing. 
Finally, let $\Omega_E$ be the least positive real period of the N\'eron 
differential of a global minimal Weierstrass equation for 
$E$, and define $C_{\infty}(E)=\Omega_E$ or $2\Omega_E$ 
according as $E(\Bbb R)$ is connected or not, and let $C_{\text{fin}}(E)$ 
denote the product of the Tamagawa factors of $E$ at the bad primes. 
The Euler product defining $L(E,s)$ converges for $\text{Re}\,s>3/2$. 
The modularity conjecture, proven by Wiles-Taylor-Diamond-Breuil-Conrad,   
implies that $L(E,s)$ has an analytic continuation to an entire function. 
The Birch and Swinnerton-Dyer conjecture relates the arithmetic data 
of $E$ to the behaviour of $L(E,s)$ at $s=1$. 

\begin{conj} (Birch and Swinnerton-Dyer) (i) $L$-function $L(E,s)$ has a 
zero of order $r=\text{rank}\;E(\Bbb Q)$ at $s=1$, 

(ii) $\sza(E)$ is finite, and 
$$ 
\lim_{s\to1}\frac{L(E,s)}{(s-1)^r} =
\frac{C_{\infty}\, (E)C_{\text{fin}}(E)\, R(E)\, |\sza(E)|}{|E(\Bbb Q)_{\text{tors}}|^2}. 
$$
\end{conj} 
If $\sza(E)$ is finite, the work of Cassels and Tate shows that its order 
must be a square. 

The first general result in the direction of this conjecture was proven for 
elliptic curves $E$ with complex multiplication by Coates and Wiles in 1976 
\cite{CW}, who showed that if $L(E,1)\not =0$, then the group $E(\Bbb Q)$ is finite. 
Gross and Zagier \cite{GZ} showed that if $L(E,s)$ has a first-order zero at $s=1$, 
then $E$ has a rational point of infinite order. Rubin \cite{Rub} proves that 
if $E$ has complex multiplication and $L(E,1)\not =0$, then $\sza(E)$ is finite. 
Let $g_E$ be the rank of $E(\Bbb Q)$ and let $r_E$ the order of the zero 
of $L(E,s)$ at $s=1$. Then Kolyvagin \cite{Kol} proved that, if $r_E\leq 1$, 
then $r_E=g_E$ and $\sza(E)$ is finite.  Very recently, Bhargava, Skinner 
and Zhang \cite{BhSkZ} proved that at least $66.48 \%$ of all elliptic curves 
over $\Bbb Q$, when ordered by height, satisfy the weak form of the Birch 
and Swinnerton-Dyer conjecture, and have finite Tate-Shafarevich group.

When $E$ has complex multiplication by the ring of integers of an imaginary 
quadratic field $K$ and $L(E,1)$ is non-zero, the $p$-part of the Birch and 
Swinnerton-Dyer conjecture has been established by Rubin \cite{Rub2} for 
all primes $p$ which do not divide the order of the group of roots of unity 
of $K$. Coates et al. \cite{cltz}, and Gonzalez-Avil\'es 
\cite{G-A} showed that there is 
a large class of explicit quadratic twists of $X_0(49)$ whose 
complex $L$-series does not vanish at $s=1$, and for which the full 
Birch and Swinnerton-Dyer conjecture is valid (covering the case $p=2$ 
when $K=\Bbb Q(\sqrt{-7})$). 
The deep results by Skinner-Urban (\cite{SkUr}, Theorem 2) 
allow, in specific cases (still assuming $L(E,1)$ is non-zero), to establish $p$-part of 
the Birch and Swinnerton-Dyer conjecture for elliptic curves  without complex multiplication 
for all odd primes $p$.

\bigskip 

This paper continues the authors previous investigations concerning 
orders of Tate-Shafarevich groups in quadratic twists of a given elliptic  
curve, and for the family of the Neumann-Setzer type elliptic curves.  
Here we present the results of our search for the 
(analytic) orders of Tate-Shafarevich groups for the cubic twists $E_m$ 
of $E: x^3+y^3=1$. These analytic orders $|\sza(E_m)|$ are the true 
ones if $|\sza(E_m)|$ are coprime to $6$ (by \cite{Rub2}). 
Our calculations extend those given by Zagier and Kramarz 
\cite{ZK} and by Watkins \cite{Wat}. 
Our main observations concern the asymptotic 
formulae in sections 3 (frequency of orders of $\sza$) and 4 
(asymptotics for the sums $\sum |\sza(E_m)|$ in the 
rank zero case), and the distributions of 
$\log L(E_m,1)$ and $\log(|\sza(E_m)|/\sqrt[t] m)$ ($t=2,3$) in sections 6 and 7. 
In section 8 we propose a variant of the asymptotic formula from section 3 
for the class numbers of real quadratic fields.

\bigskip 

This research was supported in part by PL-Grid Infrastructure.
Our computations were carried out in 2016 and 2017 on the Prometheus
supercomputer via PL-Grid infrastructure.

\section{Formula for the order of $\sza(E_m)$, when $L(E_m,1)\not=0$} 

Let $m$ be any cubefree positive integer. 
Let $E_m: x^3+y^3=m$ denote the cubic twist of $E: x^3+y^3=1$. 
Then it is plain to see that $E_m$ has the Weierstrass equation 
$y^2=x^3-432m^2$, and $E_1=X_0(27)$. 
Let $L(E_m,s)=\sum_{n=1}^{\infty}a_m(n)n^{-s}$ ($\text{Re}(s)>3/2$) denote its $L$-series. 
If $L(E_m,1)\not =0$, then the analytic order 
of $\sza(E_m)$ may be expressed as follows (see \cite{ZK})

$$
|\sza(E_m)|={L(E_m,1) \cdot T_m \over C_{fin}(E_m) \cdot C_{\infty}(E_m)}, 
$$
where

\begin{enumerate}  

\item[(i)]  $T_1=9$, $T_2=4$, and $T_m=1$ for any $m>2$; 

\item[(ii)]  $C_{\infty}(E_m)={\Gamma({1\over 3})^3\over 2\pi\sqrt{3}\sqrt[3]{m}}$  if 
$9\not | m$, and 
$C_{\infty}(E_m)={3\Gamma({1\over 3})^3\over 2\pi\sqrt{3}\sqrt[3]{m}}$  if 
$9 | m$; 

\item[(iii)]  $C_{fin}(E_m)=\prod_pC_p(E_m)$, where $C_p(E_m)=1$ if $p\not |m$, 
$C_p(E_m)=2\pm 1$ if $p\equiv \pm 1 (\text{mod} \,3)$, $C_3(E_m)=3$ if 
$m\equiv \pm 1 (\text{mod}\,9)$, $C_3(E_m)=2$ if 
$m\equiv \pm 2 (\text{mod}\,9)$, and $C_3(E_m)=1$ if $m\equiv \pm 4 (\text{mod}\,9)$ 
or $3 | m$. 

\end{enumerate}

If $\epsilon(E_m)=+1$, then the central $L$-value $L(E_m,1)$ is given by sum 
of the approximating series 
$$
L(E_m,1)=2\sum_{n=1}^{\infty}{a_m(n)\over n}e^{-2\pi n/\sqrt{N_{E_m}}},
$$
where $N_{E_m}$ is the conductor of $E_m$. The coefficients $a_m(n)$ can be computed 
as in \cite{ZK} and \cite{Wat}. In order to compute $L(E_m,1)$ with appropriate accuracy, 
we need to calculate $c\sqrt{N_{E_m}}$ terms of the approximating series (and, hence 
such a number of coefficients $a_m(n)$) for some constant $c$.

\bigskip 

\noindent 
{\it Definition}. We say, that a positive cube-free  integer $d$ satisfies condition $(*)$, 
if $\epsilon(E_d)=+1$.

\section{Frequency of orders of $\sza$}

Our data contains values of $|\sza(E_d)|$ for all positive cubic-free integers 
$d\leq 10^8$ satisfying $(*)$. Our calculations strongly suggest that for 
any positive integer $k$ there are infinitely many positive cube-free integers $d$ 
satisfying $(*)$, such that $E_d$ has rank zero and $|\sza(E_d)|=k^2$. 
Below we will state a more precise conjecture.


Let $f(X)$ denote the number of cube-free integers 
$d\leq X$, satisfying ($*$) and 
such that $|\sza(E_d)| = 1$. Let 
$g(X)$ denote the number of cube-free integers $d\leq X$, 
satisfying ($*$) and such that $L(E_d,1) = 0$. 
We obtain the following graph of the function $f(X)/g(X)$. 

\begin{figure}[H]  
\centering 
\includegraphics[trim = 0mm 20mm 0mm 15mm, clip, scale=0.4]{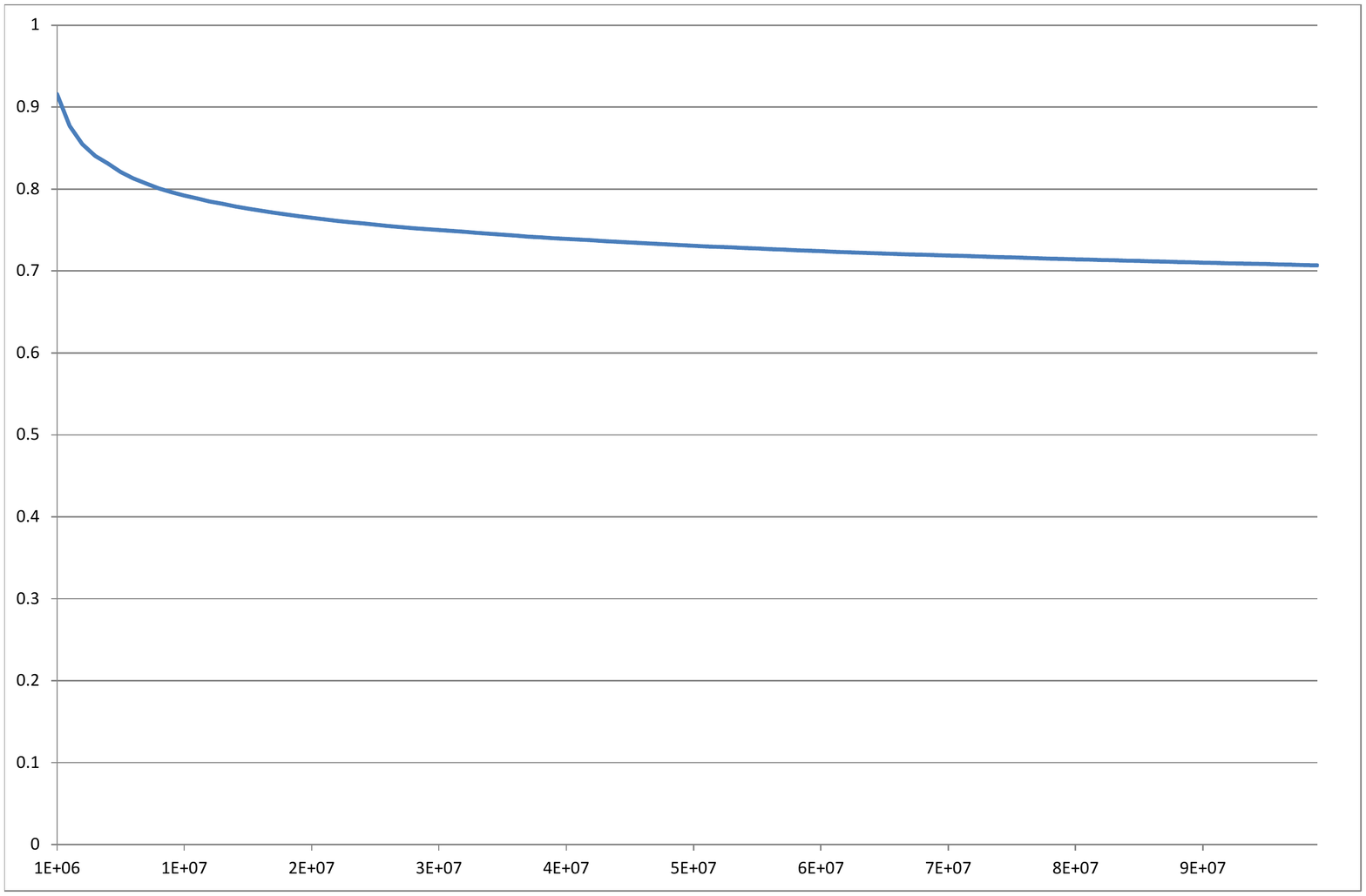}
\caption{Graph of the function $f(X)/g(X)$.} 
\end{figure}

We expect that $f(X)/g(X)$ tends to a constant ($\approx 0.7$). 
Using (\cite{Wat}, Question 1.4.1, and \cite{CKRS}), we believe the 
following asymptotic formula holds 
$$
g(X) \sim c\cdot X^{5/6}(\log X)^d, \quad X\to \infty, 
$$
with some positive $c$ and real $d$. We therefore expect a similar asymptotic 
formula for $f(X)$. Compare a similar phenomena 
for the cases of quadratic twists of elliptic curves \cite{djs} \cite{dsz2} and a family of 
Neumann-Setzer type elliptic curves \cite{dsz1}. 

\bigskip 

{\bf Remark.} Watkins claims (\cite{Wat}, Question 1.4.1 and comments after it), that if we restrict to 
the cubic twists by primes congruent to $1$ modulo $9$, then we can take 
$d=-5/8$ and $c\approx 1/6 \approx 0.16666$. 
Our calculations suggest (see the figures below) that the constant $c$ is $\approx 0.175$.  Let 
$g^*(X)$ denote the number of primes $d\leq X$, 
satisfying ($*$) and such that $L(E_d,1) = 0$. 

\begin{figure}[H]  
\centering 
\includegraphics[trim = 0mm 20mm 0mm 15mm, clip, scale=0.4]{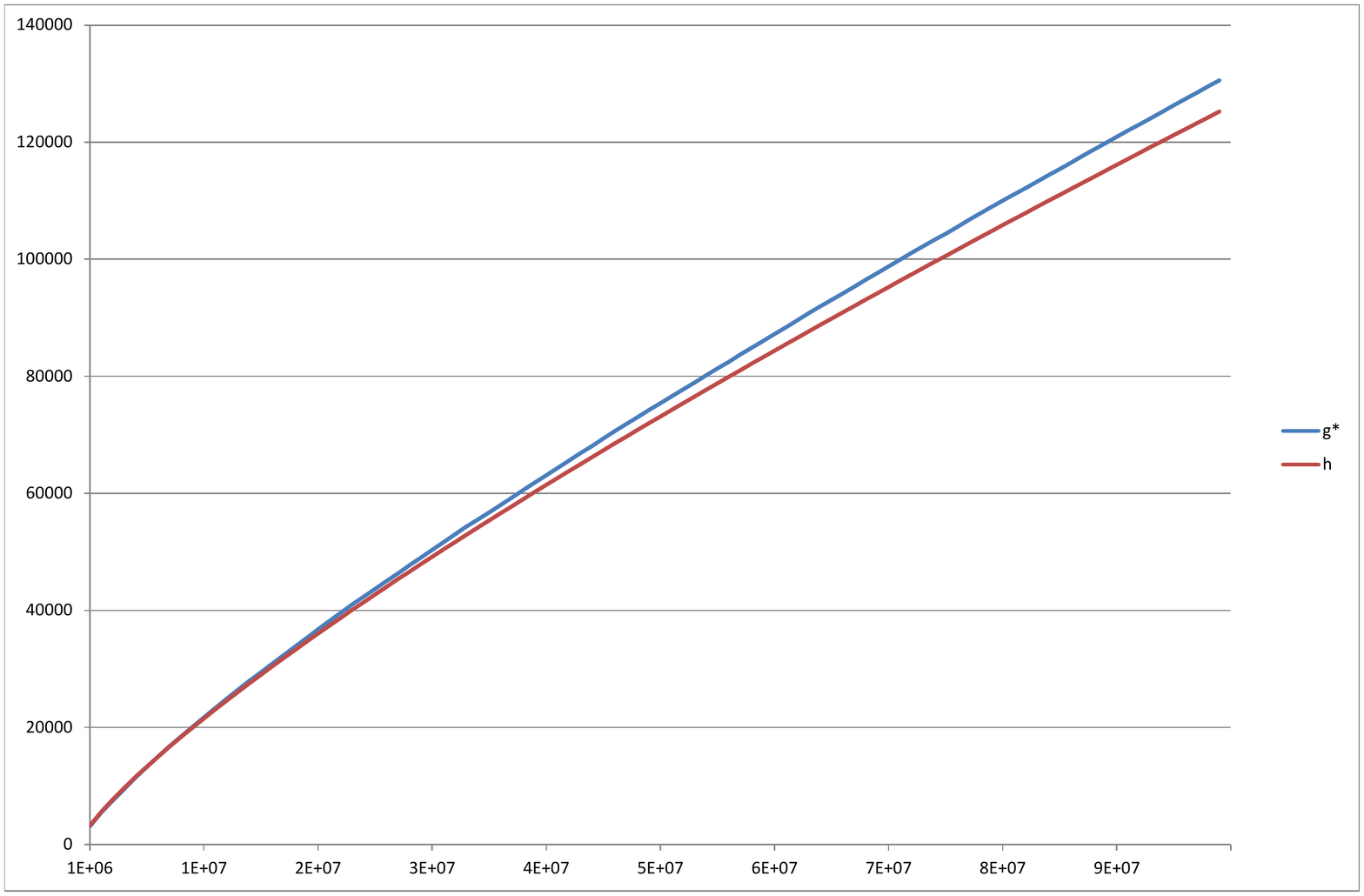} 
\caption{Graphs of the functions $g^*(X)$ and $h(X)={1\over 6}  X^{5/6}(\log X)^{-5/8}$ .} 
\end{figure}

\begin{figure}[H]  
\centering 
\includegraphics[trim = 0mm 20mm 0mm 15mm, clip, scale=0.4]{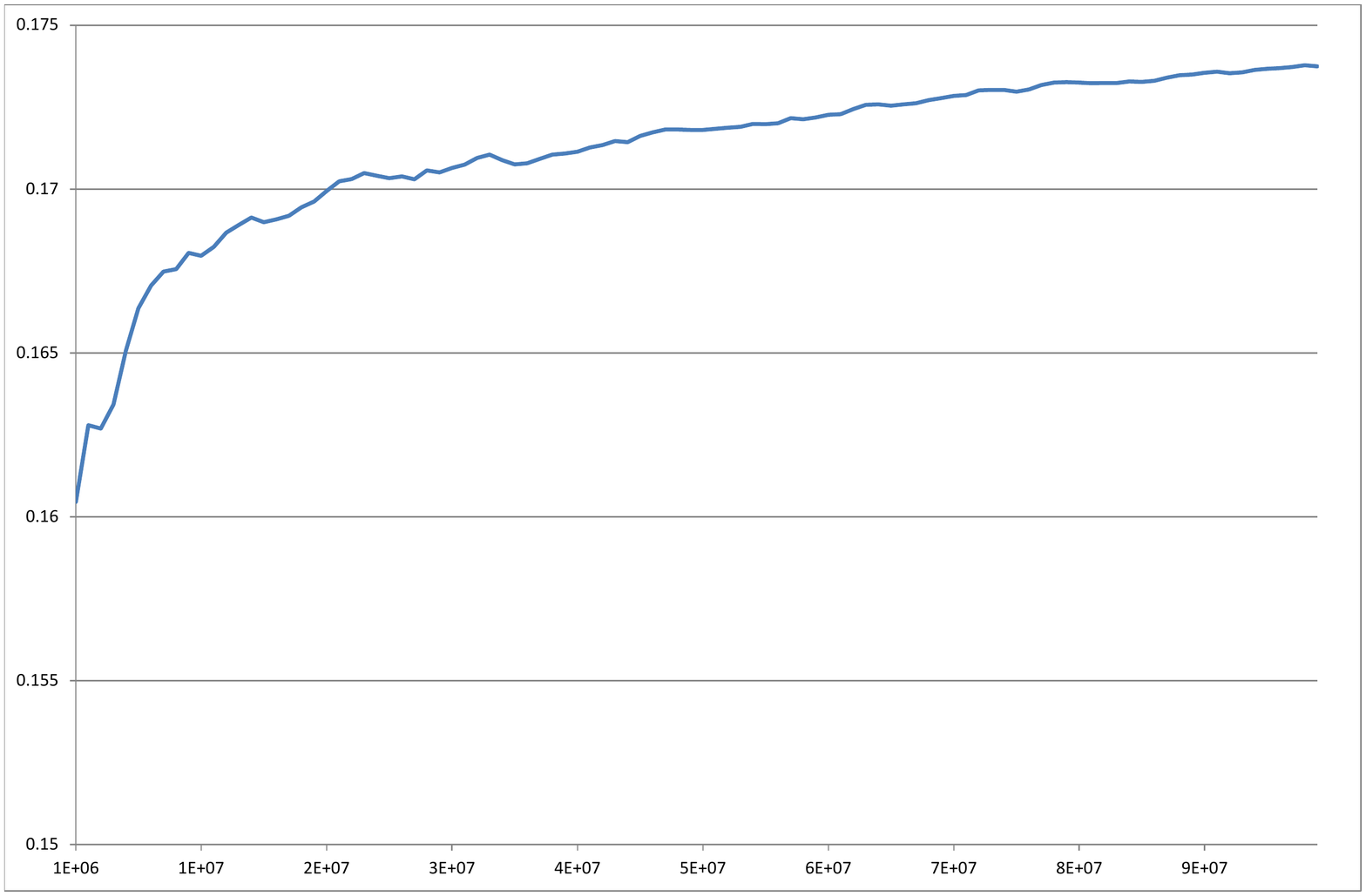} 
\caption{Graph of the function $g^*(X)/(X^{5/6}(\log X)^{-5/8})$ .} 
\end{figure}

Let us also include the graph of the function $g(X)/(X^{5/6}(\log X)^{-5/8})$ . 

\begin{figure}[H]  
\centering 
\includegraphics[trim = 0mm 20mm 0mm 15mm, clip, scale=0.4]{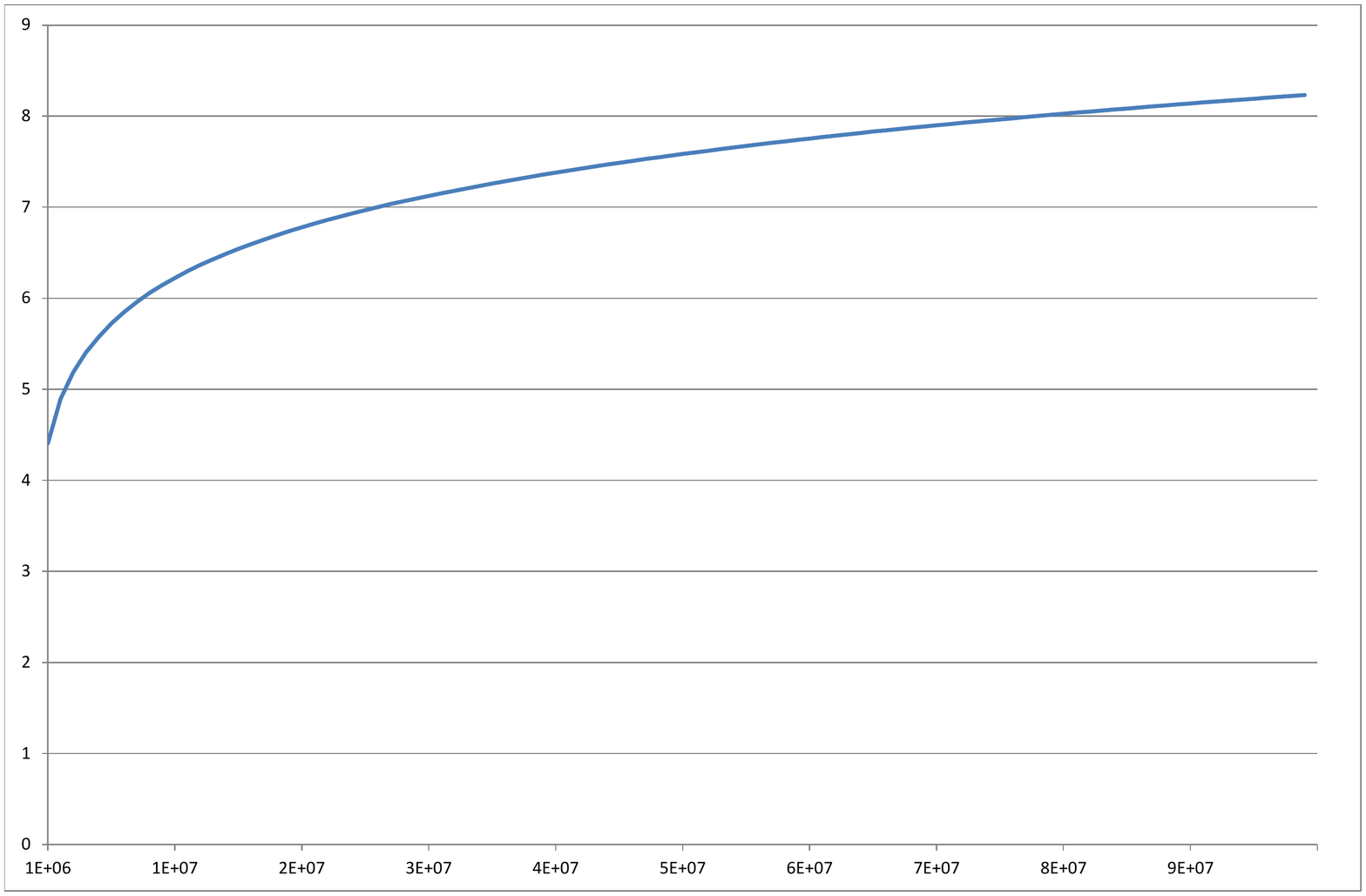} 
\caption{Graph of the function $g(X)/(X^{5/6}(\log X)^{-5/8})$ .} 
\end{figure}

\bigskip

Now let $f(k,X)$ denote the number of cube-free integers 
$d\leq X$, satisfying ($*$) and 
such that $|\sza(E_d)| = k^2$. 
Let $F(k,X):={f(X)\over f(k,X)}$. We obtain the following graphs 
of the functions $F(k,X)$ for $k=2,3,4,5,6,7$. 

\begin{figure}[H]  
\centering
\includegraphics[trim = 0mm 20mm 0mm 15mm, clip, scale=0.4]{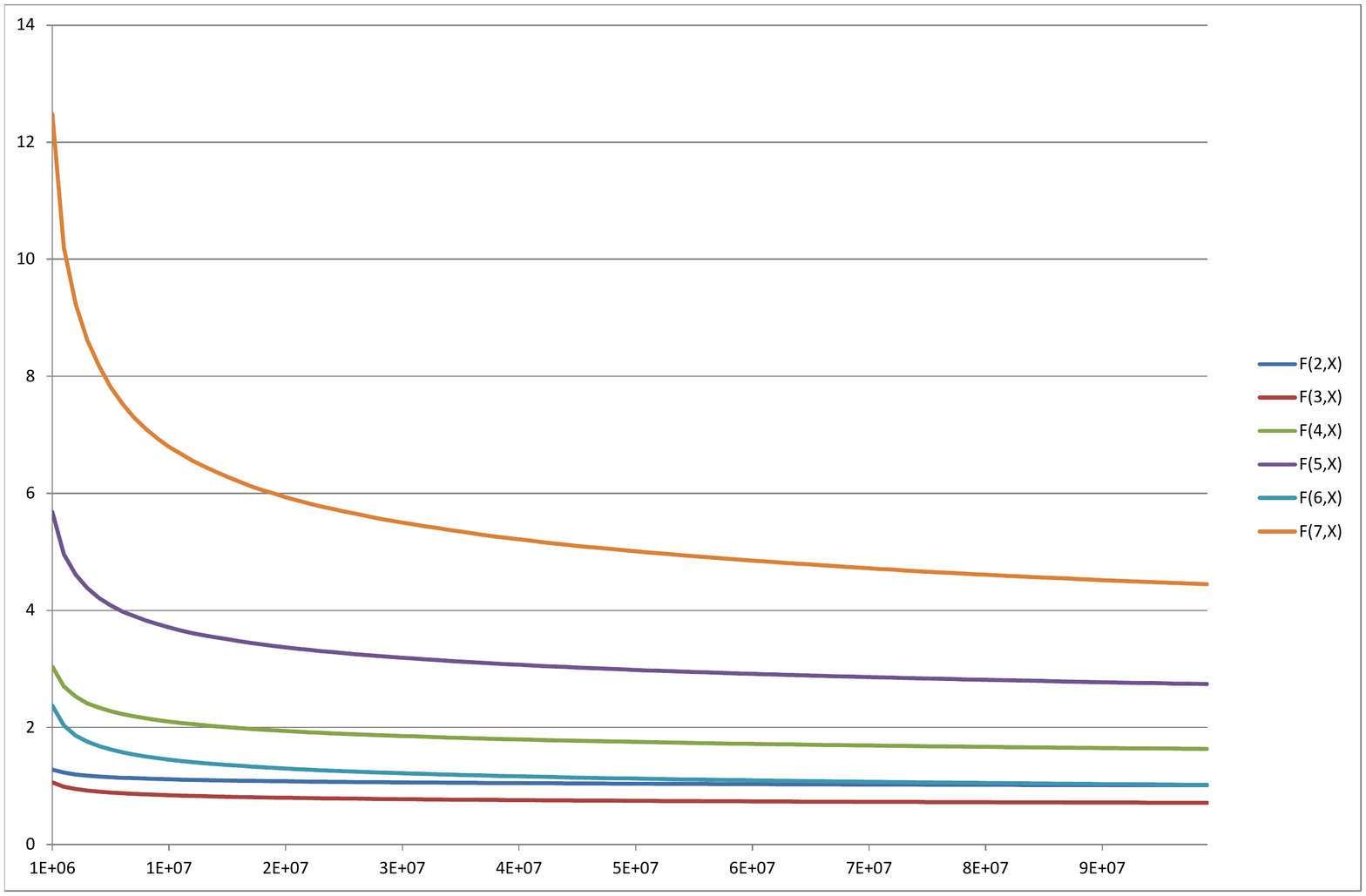}
\caption{Graphs of the functions $F(k,X)$ for $k=2,3,4,5,6,7$.} 
\end{figure}

The above calculations suggest the following general conjecture (compare 
 \cite{djs} \cite{dsz2} for the case of quadratic twists of elliptic curves, and 
\cite{dsz1} for the case of a family of Neumann-Setzer type elliptic curves). 

\bigskip 
\noindent 
{\bf Conjecture.} For any positive integer $k$ there are constants $c_k\geq 0$ anf $d_k$ 
such that 
$$
f(k,X) \sim c_k X^{5/6}(\log X)^{d_k},  \quad X\to\infty. 
$$

\bigskip 
\noindent 
{\bf Remark.} Park, Poonen, Voight and Wood \cite{PPVW}  have formulated an analogous 
(but less precise) conjecture for the family of all elliptic curves over the rationals, ordered 
by height.

\section{Variant of Delaunay's asymptotic formula}

Let 
$M^*(T):={1\over T^*}\sum |\sza(E_d)|$, 
where the sum is over primes   
$d\leq T$, satisfying (*) and $L(E_d,1)\not=0$, 
and $T^*$ denotes the number of 
terms in the sum. Similarly, let 
$N^{**}(T):={1\over T^{**}}\sum |\sza(E_d)|$, 
where the sum is over positive cube-free  integers  
$d\leq T$, satisfying (*) and $L(E_d,1)\not=0$, 
and $T^{**}$ denotes the number of terms in the sum. 
Let $f(T):={M^*(T)\over T^{1/2}}$, and 
$g(T):={N^{**}(T)\over T^{1/2}}$. We obtain the following 
picture

\begin{figure}[H]  
\centering
\includegraphics[trim = 0mm 20mm 0mm 15mm, clip, scale=0.4]{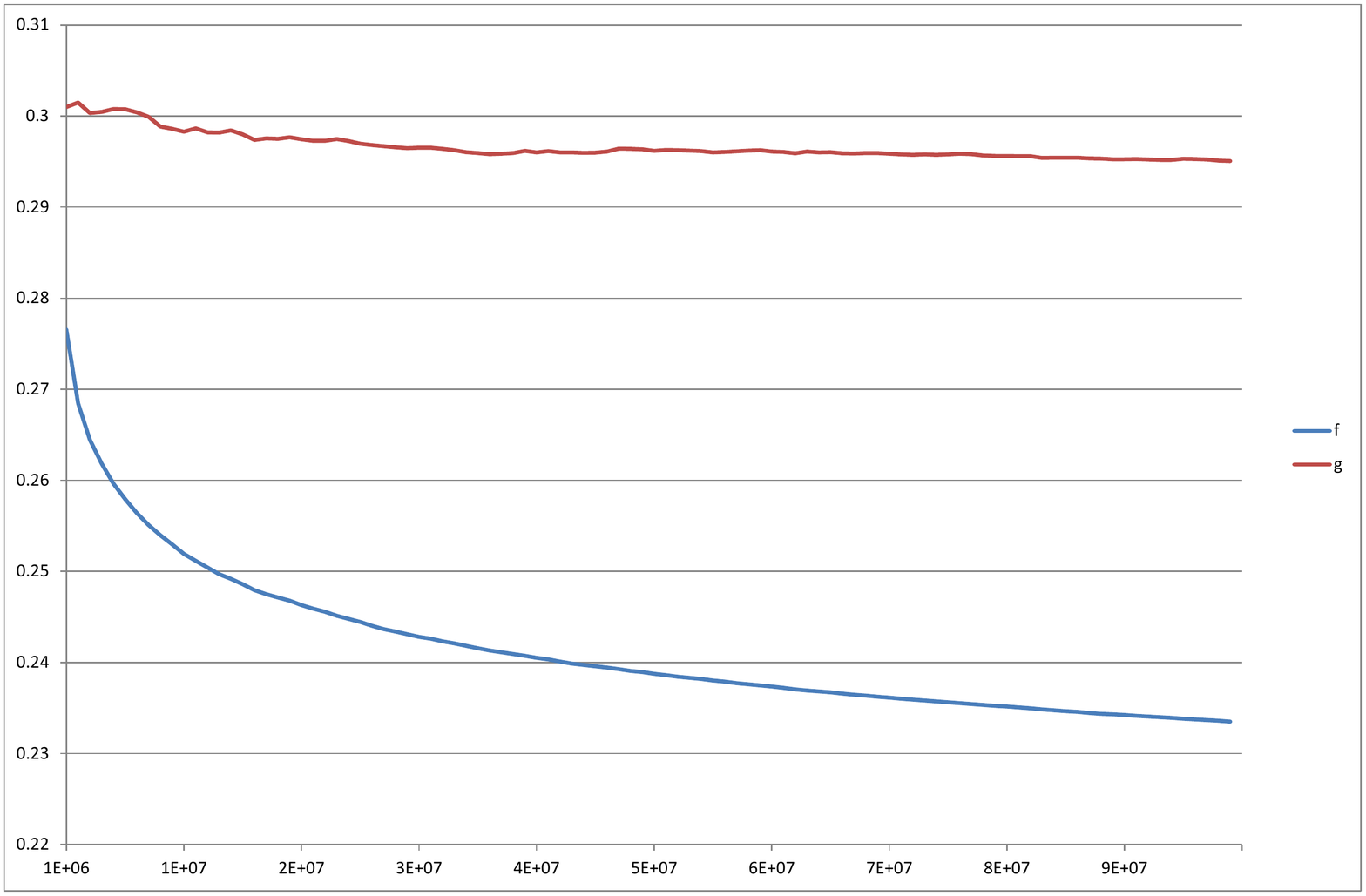} 
\caption{Graphs of the functions $f(T)$ and $g(T)$.} 
\end{figure}

Note similarity with the predictions by Delaunay \cite{Del0} 
for the case of quadra\-tic twists of a given elliptic curve (and 
numerical evidences in \cite{djs} \cite{dsz2}), and with a variant of this phenomenon 
in the case of the family of Neumann-Setzer type elliptic curves \cite{dsz1}.

\section{Cohen-Lenstra heuristics for the order of $\sza$}

Delaunay \cite{Del} have considered Cohen-Lenstra heuristics for the order of 
Tate-Shafarevich group. He predicts, among others, that in the rank zero case, 
the probability that $|\sza(E)|$ of a given elliptic curve $E$ over $\Bbb Q$ is 
divisible by a prime $p$ should be 
$
f_0(p):=1-\prod_{j=1}^{\infty}(1-p^{1-2j})={1\over p}+{1\over p^3}+... \,.    
$
Hence, $f_0(2)\approx 0.580577$, $f_0(3)\approx 0.360995$, $f_0(5)\approx 0.206660$, 
$f_0(7)\approx 0.145408$,  and so on. 

Let $F(X)$  denote the number of cube-free $d\leq X$   satisfying $(*)$ and $L(E_d,1)\not=0$, 
and let $F_p(X)$  denote the number of  such $d$'s  satisfying $p| |\sza(E_d)|$. 
Let $f_p(X):={F_p(X)\over F(X)}$, We obtain the following table (in the last row we restrict 
to prime twists)

\begin{center}
\footnotesize
\begin{longtable}{|r|r|r|r|r|r|} 
\hline 
\multicolumn{1}{|c|}{$X$} 
& 
\multicolumn{1}{|c|}{$g_2(X)$} 
& 
\multicolumn{1}{|c|}{$g_3(X)$} 
& 
\multicolumn{1}{|c|}{$g_5(X)$} 
& 
\multicolumn{1}{|c|}{$g_7(X)$} 
& 
\multicolumn{1}{|c|}{$g_{11}(X)$} 
\\ 
\hline 
\endhead 
\hline 
\multicolumn{6}{|r|}{{Continued on next page}} \\
\hline
\endfoot

\hline 
\hline
\endlastfoot 

10000000 & 0.4574860107 & 0.4528351278 & 0.0797229512 & 0.0365187357 & 0.0107055908 \\ 
20000000 & 0.4667861427 & 0.4665902606 & 0.0856954224 & 0.0406883829 & 0.0126964802 \\ 
30000000 & 0.4720389372 & 0.4743395107 & 0.0891666909 & 0.0430854869 & 0.0138608186 \\ 
40000000 & 0.4755325884 & 0.4797263355 & 0.0916462006 & 0.0448302849 & 0.0147494390 \\ 
50000000 & 0.4782835292 & 0.4838047688 & 0.0935546233 & 0.0461842060 & 0.0154253689 \\ 
60000000 & 0.4804365024 & 0.4870412651 & 0.0950607348 & 0.0472714454 & 0.0160042804 \\ 
70000000 & 0.4821166758 & 0.4897452073 & 0.0963909035 & 0.0482264317 & 0.0164998297 \\ 
80000000 & 0.4836581573 & 0.4920588749 & 0.0974999561 & 0.0490436597 & 0.0169344117 \\ 
90000000 & 0.4849849695 & 0.4940653891 & 0.0984979769 & 0.0497487127 & 0.0173190511 \\ 
100000000 & 0.4861728066 & 0.4958441463 & 0.0993871375 & 0.0503845401 & 0.0176658729 \\ 
100000000 & 0.5474977246 & 0.0713684943 & 0.1628461726 & 0.0993604813 & 0.0467913704 \\ 

\end{longtable}
\end{center}

The numerical values of $f_3(X)$ exceed the expected value $f_0(3)$ , but for $p\not=3$ 
the values $f_p(X)$ seem to tend to $f_0(p)$; additionally restricting to prime twists 
tends to speed convergence to the expected values.

\section{Distributions of $L(E_m,1)$}

It is a classical result (due to Selberg) that the values of 
$\log |\zeta({1\over 2}+it)|$ follow a normal distribution. 

Let $E$ be any elliptic curve defined over $\Bbb Q$. 
Let $\cal E$ denote the set of all fundamental discriminants $d$ with $(d,2N_E)=1$ and 
$\epsilon_E(d)=\epsilon_E\chi_d(-N_E)=1$, where $\epsilon_E$ is the root number of $E$ 
and $\chi_d=(d/\cdot)$.  Keating and Snaith \cite{KS2} have conjectured that, for 
$d\in\cal E$, the quantity $\log L(E_d,1)$ has a normal distribution with mean 
$-{1\over 2}\log \log|d|$ and variance $\log \log|d|$. 

Below we consider the case of cubic twists $E_m$ of $E=X_0(27)$. 
Our data suggest that the 
values $\log L(E_m,1)$ also follow an approximate normal distribution. 
Let $W_E = \{m \leq 10^8 : m \, \, 
\text{satisfies ($*)$} \}$ 
and $I_x = [x,x+0.1)$ for $x\in\{ -10, -9.9, -9.8, \ldots, 10 \}$.
We create histograms with bins $I_x$ from the data 
$\left\{ \left(\log L(E_m,1) + \frac 1 2 \log\log m\right)/\sqrt{\log\log m} : m\in W_E\right\}$. 
Below we picture this histogram.

\begin{figure}[H]  
\centering
\includegraphics[trim = 0mm 20mm 0mm 15mm, clip, scale=0.4]{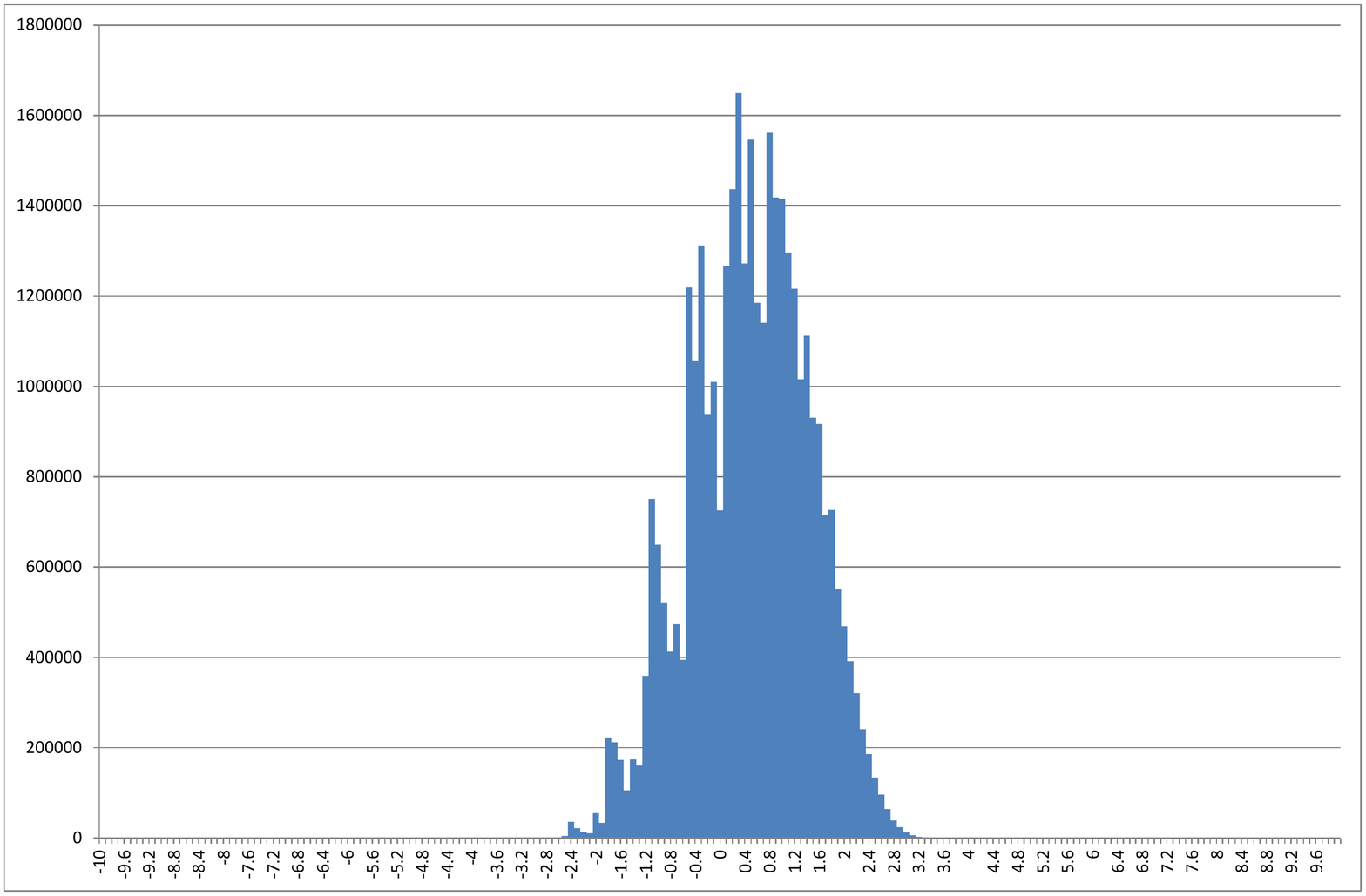} 
\caption{Histogram of values $\left( \log L(E_m,1) + {1\over 2} 
\log\log m \right)/\sqrt{\log\log m}$ for $m\in W_E$.} 
\end{figure}

\section{Distributions of $|\sza(E_m)|$}

It is an interesting question to find 
results (or at least a conjecture) on distribution of the order of the Tate-Shafarevich 
group in family of elliptic curves. 
It turns out that in a case of rank zero quadratic twists $E_d$ of a fixed elliptic curve $E$ 
the values of $\log(|\sza(E_d)|/\sqrt{d})$ are the natural ones to consider (compare  
the numerical experiments in \cite{djs}, \cite{dsz2}). We also have good conjecture 
for a family of rank zero Neumann-Setzer type elliptic curves \cite{dsz1}. 

Let us consider the family $E_m$ of cubic twists of $E=X_0(27)$. In this case we will 
create histograms for the values  $\log(|\sza(E_m)|/\sqrt[t] m)$, $t=2, 3$ separately.  
Let $W_E = \{m \leq 10^8 : m \, \, 
\text{satisfies ($*)$} \}$ 
and $I_x = [x,x+0.1)$ for $x\in\{ -10, -9.9, -9.8, \ldots, 10 \}$.
We create histograms with bins $I_x$ from the data 
$\left\{ \left(\log(|\sza(E_m)|/\sqrt[t] m)  + \frac 1 2 \log\log m\right)/\sqrt{\log\log m} : m\in W_E\right\}$.  
 Below we picture these histograms.

 \begin{figure}[H]  
\centering
\includegraphics[trim = 0mm 20mm 0mm 15mm, clip, scale=0.4]{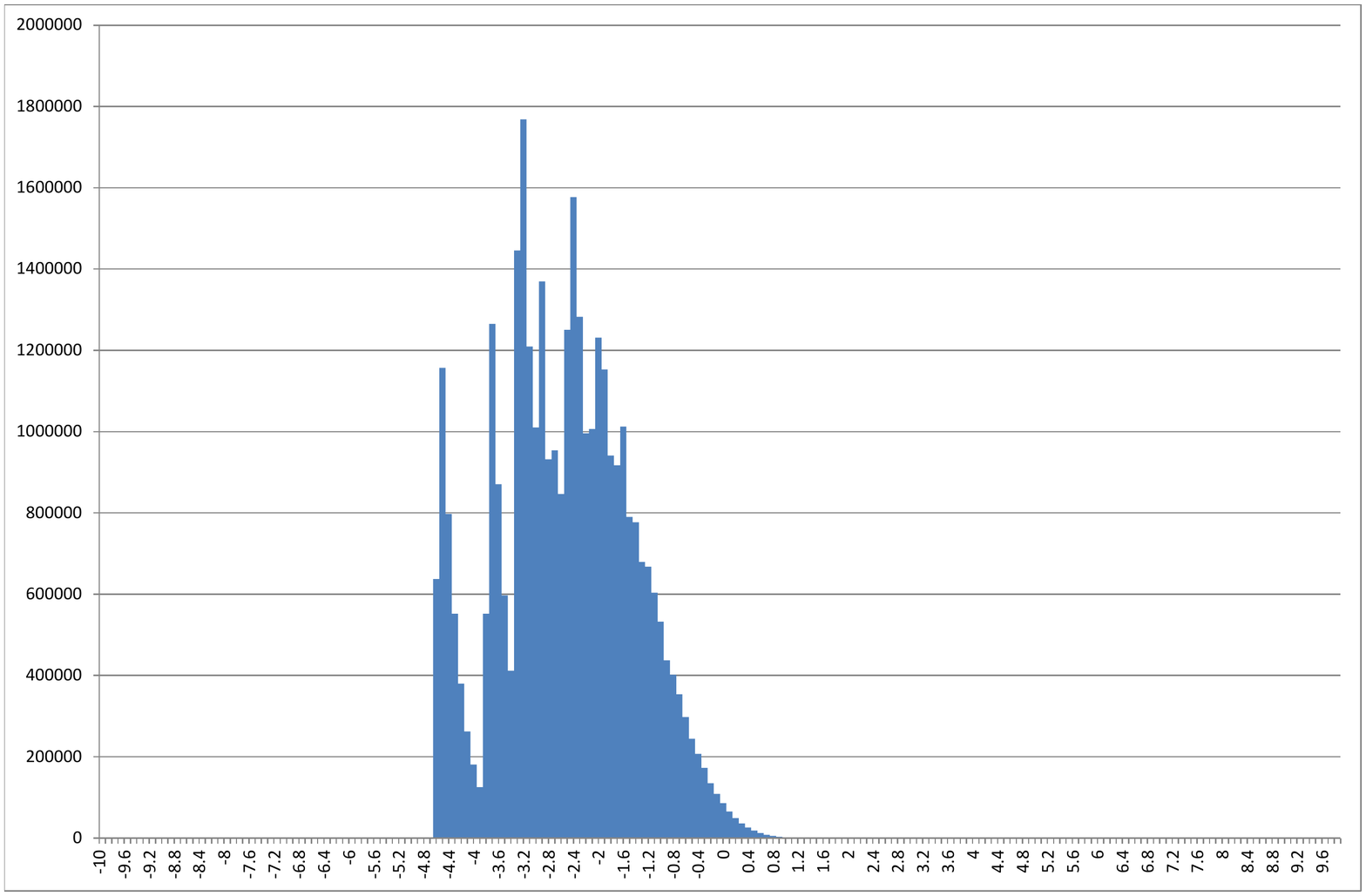} 
\caption{Histogram of values $\left( \log(|\sza(E_m)|/\sqrt m)  +\frac 1 2 \log\log m \right)/\sqrt{\log\log m}$ 
for $m\in W_E$.} 
\end{figure}

\begin{figure}[H]  
\centering
\includegraphics[trim = 0mm 20mm 0mm 15mm, clip, scale=0.4]{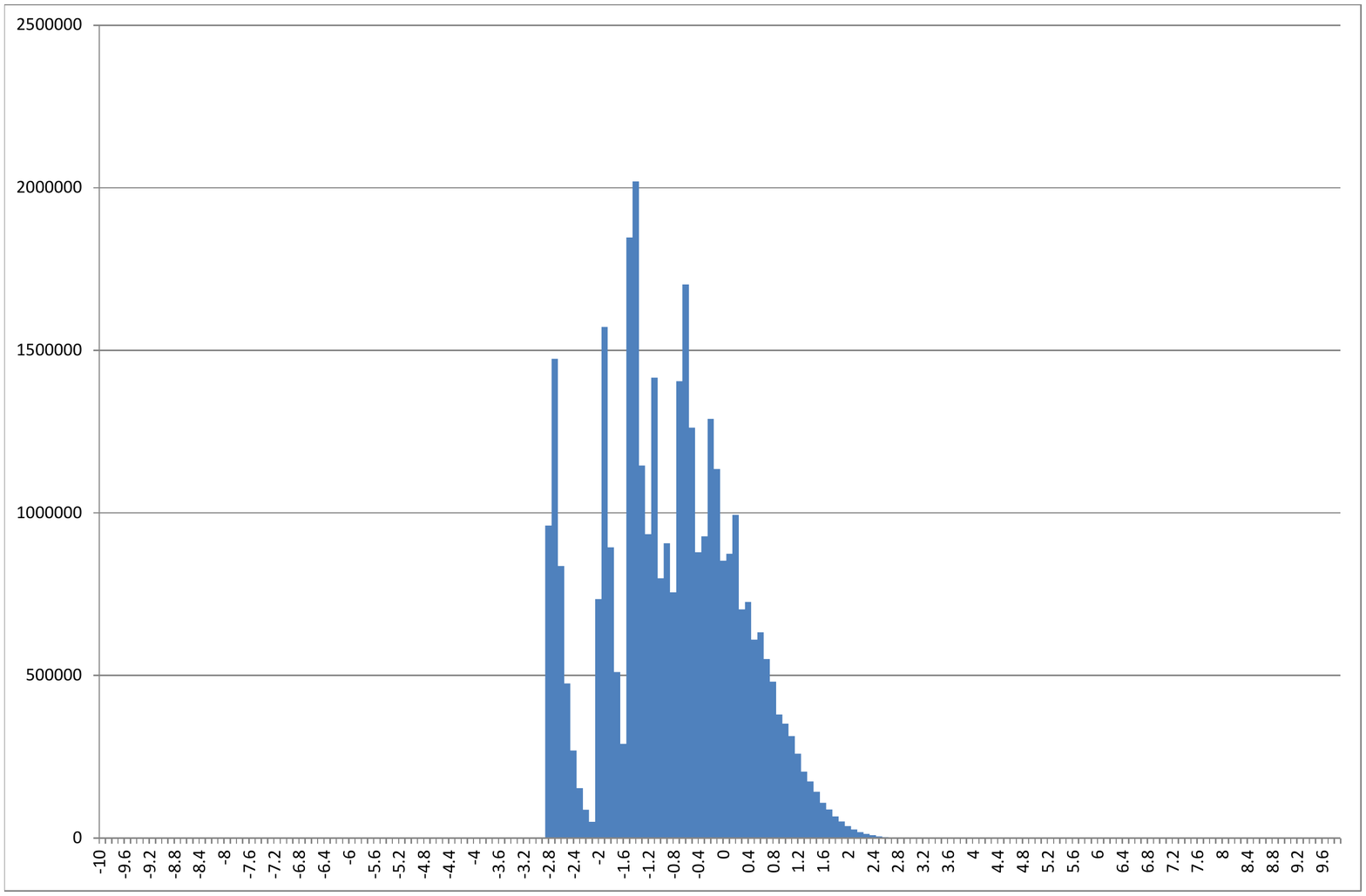} 
\caption{Histogram of values $\left( \log(|\sza(E_m)|/\sqrt[3] m)  + \frac 1 2 \log\log m \right)/\sqrt{\log\log m}$ 
for $m\in W_E$.} 
\end{figure}

\section{Observations concerning the class numbers of real quadratic fields}

Consider a real quadratic field $K=\Bbb Q(\sqrt{d})$ ($d$ a positive square-free integer); let $h(d)$ 
denote its class number. We calculated the values $h(d)$ for all positive square-free integers 
$d \leq 3\cdot 10^{10}$. Our observations suggest that $h(d)$'s behave in a similar way to the 
orders of Tate-Shafarevich groups in some families of rank zero elliptic curves (i.e. quadratic or cubic twists 
of a given one).  

Let $h(k,X)$ denote the number of positive square-free integers  
$0<d\leq X$  such that $h(d) = k$. 
Let $H(k,X):={h(1,X)\over h(k,X)}$. We obtain the following graphs 
of the functions $H(k,X)$ for $2\leq k\leq 10$.

\begin{figure}[H]  
\centering
\includegraphics[trim = 0mm 20mm 0mm 15mm, clip, scale=0.4]{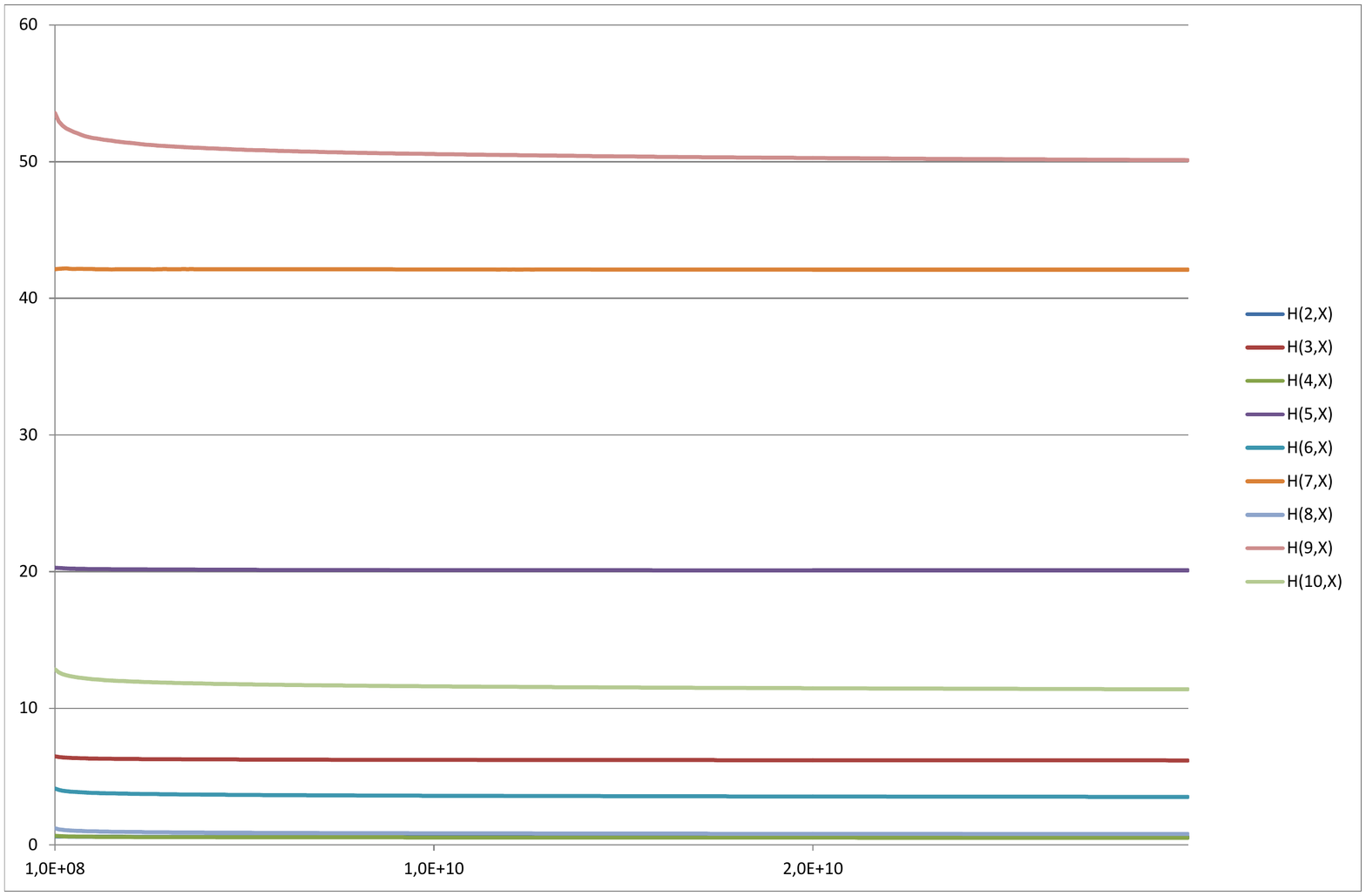} 
\caption{Graphs of the functions $H(k,X)$ for $2\leq k\leq 10$.} 
\end{figure}

Now let us consider graphs of the functions $h(1,X)/(X^{5/6}(\log  X)^r)$, $r=0,1$. 

\begin{figure}[H]  
\centering
\includegraphics[trim = 0mm 20mm 0mm 15mm, clip, scale=0.4]{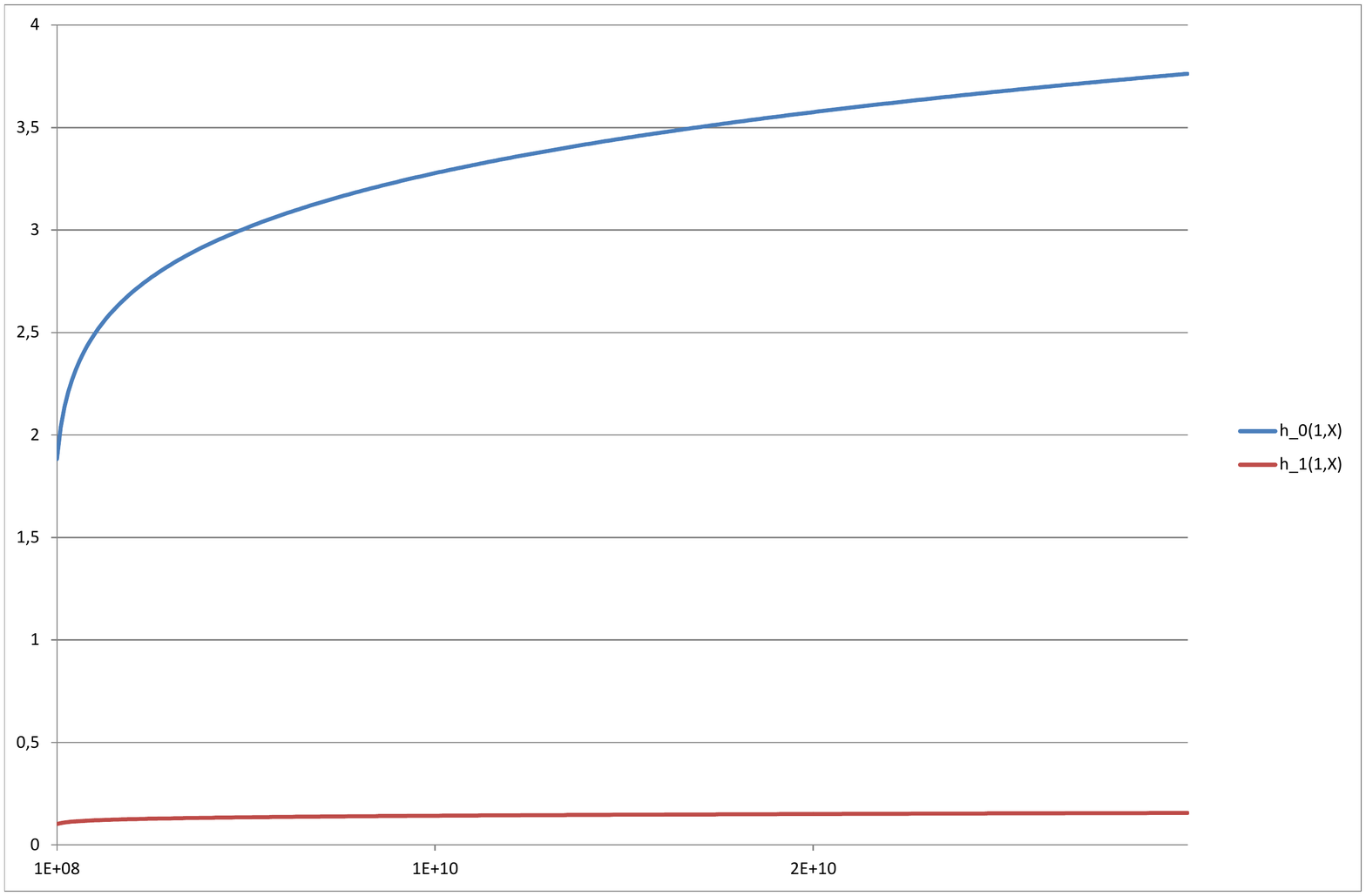} 
\caption{Graphs of the functions $h_r(1,X):=h(1,X)/(X^{5/6}(\log  X)^r)$, $r=0,1$.} 
\end{figure}

The above calculations suggest the following (optimistic) conjecture. 

\bigskip 
\noindent 
{\bf Conjecture.} For any positive integer $k$ there are positive constants $r_k$, $s_k$ such that 
$$
h(k,X) \sim r_kX^{5/6}(\log X)^{s_k},  \quad X\to\infty. 
$$

\bigskip 
\noindent 
{\bf Remark.} The Gauss' class-number one problem for real quadratic fields states that there are infinitely 
real quadratic fields with trivial ideal class group. It is still an open problem; note that it is not even 
known if there are infinitely many number fields with a given class number. Therefore the above conjecture 
is a highly optimistic version of these open  questions.

Institute of Mathematics, University of Szczecin, Wielkopolska 15, 
70-451 Szczecin, Poland; E-mail addresses: andrzej.dabrowski@usz.edu.pl and dabrowskiandrzej7@gmail.com;  
lucjansz@gmail.com


\begin{thebibliography}{999} 

\bibitem{BhSkZ} M. Bhargava, Ch. Skinner, W. Zhang, {\it A majority of elliptic 
curves over $\mathbb Q$ satisfy the Birch and Swinnerton-Dyer conjecture}, 
arxiv.org/abs/1407.1826 

\bibitem{cltz} J. Coates, Y. Li, Y. Tian, S. Zhai, {\it Quadratic 
twists of elliptic curves}, Proc. London Math. Soc. {\bf 110} (2015), 357-394  

\bibitem{CW} J. Coates, A. Wiles, {\it On the conjecture of Birch and 
Swinnerton-Dyer}, Invent. Math. {\bf 39} (1977), 223-251 

\bibitem{CKRS} J. B. Conrey, J. P. Keating, M. O. Rubinstein, N. C. Snaith, 
{\it On the frequency of vanishing of quadratic twists of modular $L$-functions}. 
In: Number Theory for the millennium, I, (Urbana, IL, 2000) (ed. by M. A. Bennett et al.), 301-315 


\bibitem{djs} A. D\k{a}browski, T. J\k{e}drzejak, L. Szymaszkiewicz, 
{\it Behaviour of the order of Tate-Shafarevich groups for the 
quadratic twists of $X_0(49)$}, In: Elliptic Curves, Modular Forms 
and Iwasawa Theory (in honour of John Coates' 70th birthday), Springer 
Proceedings in Mathematics and Statistics {\bf 188} (2016), 125-158  

\bibitem{dsz1} A. D\k{a}browski, L. Szymaszkiewicz, {\it Orders of Tate-Shafarevich 
groups for the Neumann-Setzer type elliptic curves},   Math. Comput. {\bf 87} (2018) 

\bibitem{dsz2} A. D\k{a}browski, L. Szymaszkiewicz, {\it Behaviour of the order 
of Tate-Shafarevich groups for the quadratic twists of elliptic curves},  arXiv:1611.07840 


\bibitem{Del0} C. Delaunay, {\it Moments of the orders of Tate-Shafarevich 
groups}, Int. J. Number Theory {\bf 1} (2005), 243-264  

\bibitem{Del} C. Delaunay, {\it Heuristics on class groups and on Tate-Shafarevich 
groups: the magic of the Cohen-Lenstra heuristics}. In: Ranks of elliptic curves 
and random matrix theory, London Math. Soc. Lecture Ser. {\bf 341} (2007), 
323-340 

\bibitem{G-A} C.D. Gonzalez-Avil\'es, {\it On the conjecture of Birch 
and Swinnerton-Dyer}, Trans. Amer. Math. Soc. {\bf 349} (1997), 4181-4200 

\bibitem{GZ} B. Gross, D. Zagier, {\it Heegner points and derivatives 
of $L$-series}, Invent. Math. {\bf 84} (1986), 225-320 



\bibitem{KS2} J.P. Keating, N.C. Snaith, {\it Random matrix theory and $\zeta(1/2+it)$}, 
Comm. Math. Phys. {\bf 214}(1) (2000), 57-89 

\bibitem{Kol} V. Kolyvagin, {\it Finiteness of $E(\Bbb Q)$ and $\sza(E)$ 
for a class of Weil curves}, Math. USSR Izv. {\bf 32} (1989), 523-541 




\bibitem{PARI} 
The PARI~Group, PARI/GP version \texttt{2.7.2}, Bordeaux, 2014, \texttt{http://pari.math.u-bordeaux.fr/} 

\bibitem{PPVW} J. Park, B. Poonen, J. Voight, M. M. Wood, {\it A heuristic for boundedness of 
ranks of elliptic curves}, www-math.mit.edu/~poonen/papers/bounded-ranks.pdf 


\bibitem{Rub} K. Rubin, {\it Tate-Shafarevich groups and $L$-functions 
of elliptic curves with complex multiplication}, Invent. Math. {\bf 89} 
(1987), 527-560 

\bibitem{Rub2} K. Rubin, {\it The "main conjectures" of Iwasawa theory 
for imaginary quadratic fields}, Invent. Math. {\bf 103} (1991), 25-68 

\bibitem{SkUr} Ch. Skinner, E. Urban, {\it The Iwasawa main conjectures 
for $GL_2$}, Invent. Math. {\bf 195} (2014), 1-277  

\bibitem{Wat} M. Watkins, {\it Rank distributions in a family of cubic twists}. In:  
Ranks of elliptic curves and random matrix theory, 237-246, London Math. Soc. Lecture Note Ser., 
{\bf 341}, Cambridge Univ. Press, Cambridge, 2007  

\bibitem{ZK} D. Zagier, G. Kramarz, {\it Numerical investigations related to the 
$L$-series of certain elliptic curves},  J. Indian Math. Soc. (N.S.) {\bf 52} (1987), 
51-69  

\end{thebibliography}
\end{document}